%% file: 567.tex
 \documentstyle[12pt]{article}
\def\indd #1#2#3{#1 \mathbin{\mathop{\down}^d_{#2}} #3}
\def\bar#1{\overline{#1}}

\def\th{\rm Th }

\def\indd #1#2#3{#1 \mathbin{\mathop{\down}^d_{#2}} #3}
\def\bar#1{\overline{#1}}

\def\bK{{\bf K}}

\def\bbK{{\bf {\overline { K}}}_0}
\def\th{\rm Th }
\def\icl{\rm icl}

\input jbmac

\title{DOP and FCP in Generic Structures}
\author{John T. Baldwin
\thanks{Partially supported by NSF grant
9308768.}
\\ Department of Mathematics, Statistics and  Computer
Science\\University of Illinois at Chicago
\and
Saharon Shelah\thanks{This is paper 567 in Shelah's bibliography.   
Both authors thank
Rutgers University and the Binational Science Foundation for partial
support of this research.}\\ Department of Mathematics\\
Hebrew University of Jerusalem\\Rutgers University}

\begin{document}
\maketitle
\def\subi{\leq_{i}}
\def\subm{\leq_{s}}
\def\bT{{\bf T}}
\def\bK{{\bf K}}
\section{Context}
%\jtbnumpar{Context}
We work throughout in a finite relational language $L$.  This paper  
is
built on \cite{BaldwinShelahran} and \cite{BaldwinShiJapan}.  We  
repeat
some of the basic notions and results from these papers for the
convenience of the reader but familiarity with the setup in the
first few sections of \cite{BaldwinShiJapan} is needed to read this  
paper.
 Spencer and Shelah \cite{ShelahSpencer} constructed for each  
irrational $\alpha$
 between $0$ and
$1$
the theory $T^{\alpha}$ as the almost sure theory of random graphs  
with edge probability $n^{-\alpha}$.
In \cite{BaldwinShelahran} we
proved that this was the same theory as
the theory $T_{\alpha}$ built by constructing a generic model in
\cite{BaldwinShiJapan}.  In this paper we explore some of the
more subtle model theoretic properties of this theory.  We show that
$T^{\alpha}$ has the dimensional order property and does not have the  
finite cover property.

We
work in the framework of \cite{BaldwinShiJapan} so probability  
theory
is not needed in this paper. This choice allows us to consider a  
wider
class of theories than just the $T_{\alpha}$.  The basic facts cited  
from \cite{BaldwinShiJapan} were due to Hrushovski  
\cite{Hrustableplane}; a full bibliography is in  
\cite{BaldwinShiJapan}.  For general background in stability theory  
see \cite{Baldwinbook} or \cite{Shelahbook2nd}.

We  work at three levels of generality.
The first is given by an
axiomatic framework in Context~\ref{finalsetup}.
Section~\ref{indord} is carried out in this generality.
The  main family of examples for this context is described
Examples~\ref{mainexamp}.
Sections~\ref{findtype} and \ref{fcp} depend on a
function $\delta$ assigning a real number to each finite  
$L$-structure
as in these examples.
Some of the constructions in Section~\ref{findtype} (labeled at the
time) use heavily
the restriction of the  class of examples to
graphs.  The first author acknowledges useful discussions on this  
paper with Sergei Starchenko. 

\jtbnot \label{setup} Let ${\bK_0}$ be a class of
finite structures closed under substructure and
isomorphism and  containing the empty structure.  Let $\bbK$ be the
universal class determined by $\bK_0$.  % and let
% {\bK} be an arbitrary subclass of $\bbK$ which contains ${\bK_0}$.

\jtbnot
Let $B \inter C = A$.
 The {\em free amalgam} of $B$ and $C$ over $A$, denoted
$B \bigotimes_{A} C$,
is the structure with universe $BC$
but no relations not in $B$ or $C$.

We write $A \subseteq_{\omega} B$ to mean $A$ is a finite subset of  
$B$.
A structure $A$ is called {\em discrete} if there are no relations  
among the
elements of $A$.
Let $\delta: {\bf
K}_0 \mapsto \Re^+$ (the nonnegative reals) be an arbitrary function
with $\delta(\emptyset) = 0$.
Extend $\delta$ to
$d: {\bbK} \times {\bK}_0
 \mapsto \Re^+$ by for each $N \in {\bbK}$,
$$d(N,A) =
\inf \{\delta(B): A \subseteq B \subseteq_{\omega} N\}.$$  We  
usually
write $d(N,A)$ as $d_N(A)$.  We only use this definition when  
$\delta$
is defined on every finite subset of $N$.
We  will omit the subscript $N$ if it is clear from context.

For $g = \delta$ or $d_N$ and finite $A,B$, we define {\em relative  
dimension } by
$g(A/B) = g(AB) - g(B)$.  For infinite $B$ and finite $A$, $d(A/B) =  
\inf\{d(A/B_0): B_0 \subset_{\omega} B\}$.  This definition is  
justified in e.g. Section 3 of
 \cite{BaldwinShiJapan}.  For any finite sequence $\abar \in N
$, $d_N(\abar)$ is the same as $d_N(A)$ where $\abar$ enumerates  
$A$.

Consider a finite  structure $B$ for a finite
relational language $L$.  
We assume that each relation of $L$ holds of a tuple $\abar$ only if  
the elements $\abar$ are
distinct and if $R(\abar)$ holds, $R(\abar')$ holds for any  
permutation $\abar'$ of $\abar$.

$R(B)$ denotes the collection of subsets $B_0 =\{b_1, \ldots b_n\}$  
of
$B$ such that for some (any) ordering $\bbar$ of $B_0$,
$B\sat R(\bbar)$ for some relation symbol $R$ of
$L$;  $e(B)= |R(B)|$.
Let $A$, $B$, $C$ be disjoint sets.  
We
write $R(A,B)$ for the collection of subsets from $AB$
that satisfy some relation of $L$ (counting with multiplicity if a  
set satisfies
more than one relation)
  and contain at least one
member of $A$ and one of $B$.
Write $e(A,B)$ for $|R(A,B)|$.
Similarly, we
write $R(A,B,C)$ for the collection of subsets from $ABC$
that satisfy some relation of $L$ and contain at least one
member of $A$ and one of $C$.
Write $e(A,B,C)$ for $|R(A,B,C)|$.

\jtbnumpar{Example}
\label{mainexamp}
The most important examples arise by defining $\delta$ as follows.   
In the last section of \cite{BaldwinShiJapan} we enumerated several  
other examples to
which this axiomatization applies.  Let
$$\delta_{\beta,\alpha}(A) = \beta |A| - \alpha e(A).$$
We may write $\delta_{\alpha}$ for $\delta_{1,\alpha}$.
The class $\bK_{\alpha}$ is the collection of finite $L$-structures  
$A$
such that
for any $A' \subseteq A$, $\delta_{\alpha}(A') \geq 0$.
We denote by $T_{\alpha}$ the theory of the generic model of
$\bK_{\alpha}$.

\jtbnumpar{Axioms}
\label{yax} Let $N$ be in $ {\bbK}$ and let $A,\, B,\, C \in {\bf  
K_0}$
be
  substructures of $N$.

\begin{enumerate}
\item
If $A$, $B$, and $C$ are disjoint then
$\delta(C/A) \geq \delta(C/AB)$.

\item
For every $n$ there is an $\epsilon_n>0$ such that if $|A|<n$ and
$\delta(A/B) < 0$ then $\delta(A/B) \leq -\epsilon_n$.
\item  There is a real number $\epsilon$ independent of $N,A,B,C$  
such
that if $A,B,C$ are disjoint subsets of a model $N$ and
$\delta(A/B) - \delta(A/BC) < \epsilon$ then
$R(A,B,C) = \emptyset$ and $\delta(A/B) = \delta(A/BC)$.
\item For each $A \in \bK_0$, and each $A'\subseteq A$, $\delta(A')  
\geq
0$.
 \end{enumerate}

We call a function $d = d_N$ derived from $\delta$ satisfying  
Axioms~\ref{yax}
a {\em dimension function}.
\begin{lemm}
\label{yaxtrue}
% Consider the situation described in
%Definition~\ref{finalsetup}.
If $\delta$ is a dimension function
satisfying the properties of Axiom~\ref{yax} and $\subm$ ( read {\em  
strong submodel} ) is defined
by $A \subm N$ if $d_N(A) = d_A(A)$,
then $\subm$  satisfies the following
propositions.  Let $M,N,N'\in \bbK$.

{ \bf A1}.  $M\subm M$.

{ \bf A2}. If $M\subm N$ then $M\subseteq N$.

{ \bf A3}.
\(M\subm N'\subm N\) implies \( M\subm N'\).

{ \bf A4}.
If $M\subm N$, $N'\subseteq N$
%and $A\bigcap C$ is finite,
then $M\bigcap N' \subm N$.

{\bf { A5}.} For all $M\in \bbK$, $\emptyset \subm M$.
\end{lemm}
%\begin{prop}\label{deltadisjoint}
%  Let $N\in {\bbK}$, and $A,\, B\in {\bf K_0}$ be disjoint
%  substructures of $N$. If $\delta = \delta_{\alpha,\beta}$ and
%$A \inter B = \emptyset$, then
%$\delta(A/B)
% = \delta(A) - \alpha e(A,B)$.
%\end{prop}

We need to analyze extensions which are far from being strong.
\jtbdef
\label{defins1}
For $A,B \in S(\bK_0)$, $A \subi B$ if $A \subseteq B$ but there is  
no
$A'$ properly contained in $B$ with $A \subseteq A' \subm B$.  If $A  
\subi B$, we say $B$ is an
{\em intrinsic} extension of $A$.

\jtbdef The intrinsic closure of $A$ in $M$, $\icl_M(A)$ is the union  
of
$B$ with $A \subseteq B \subseteq M$ and $A \subi B$.  When $M$ is  
clear from
context, we write $\Abar$ for $\icl_M(A)$.
The intrinsic closure can be more finely analyzed as follows.
\begin{enumerate}
\item For any $M \in \bK$, any $m\in \omega$, and any $A \subseteq  
M$,
  $$\textstyle{\icl^m_M(A)} = \bigunion\{B:A \subi B \subseteq M \&
  |B-A| < m\}.$$
\item $$\textstyle{\icl_M(A)} = \bigunion_m
\icl^m_M(A).$$
\item $M$ has {\em finite closures} if for each finite $A \subseteq  
M$,
$\icl(A)$ is finite.
\item $\bK$ has {\em finite closures} if each $M \in \bK$ has finite  
closures.
\end{enumerate}

Using {\bf A4}, note that the intrinsic closure of $A$ in $M$ is 
the intersection of the strong substructures of $M$ which contain  
$A$.  Thus, when finite,
$\icl_M(A) \in \bK_0$ and is a strong substructure of $M$.   
Moreover,
a countable $M$ has finite closures if and only if $M$ can be written  
as an increasing
union of finite strong substructures.

\jtbdef
The countable
model $M\in \bbK$ is {\em $({\bK_0}, \subm)$-generic} if
\begin{enumerate}
\item
     If $A\subm M, A\subm B\in \bK_{0}$, then there exists $B'\subm M  
$
such that $ B\cong_{A} B'$,
\item
$M$ has finite closures.
\end{enumerate}

\jtbnumpar{Fact}  If $(\bK_0,\subm)$ satisfies the properties of
Lemma~\ref{yaxtrue} and the amalgamation property with respect to
$\subm$ then there is a countable $\bK_0$-generic model.

 \jtbnumpar{Context}
\label{finalsetup}
Henceforth, $(\bK_0,\subm)$ is class of finite structures closed  
under
isomorphism and substructure with $\subm$ induced by a function  
$\delta$
obeying Axioms~\ref{yax}.  Moreover, we assume $(\bK_0,\subm)$  
satisfies
the amalgamation property and $\bK$ is the class of models of the  
theory
of the generic model $M$
of $(\bK_0,\subm)$.  $\Mscr$ is a large saturated
model of $T = \th(M)$.  In the absence of other specification,
the dimension function $d$ is the  function
induced on $\Mscr$ by $\delta$ and we work with substructures of  
$\Mscr$.

\section{Independence and Orthogonality}
\label{indord}

As indicated in Context~\ref{finalsetup}, the following definitions  
take place
in a suitably saturated model elementarily equivalent to the generic.   
We work in that context throughout this section.

\jtbdef
We say the finite sets $A$ and $B$ are $d$-independent over $C$
and write\begin{enumerate}
\item
$\indd A C B$ if
\begin{enumerate}
\item $d(A/C) = d(A/CB)$.
\item $\bar{AC} \inter \bar{BC} \subseteq \Cbar$.
\end{enumerate}
\item
We say the (arbitrary) sets $A$ and $B$ are $d$-independent over $C$
and write $\indd A C B$ if for every finite $A' \subseteq A$ and
$B'\subseteq B$, $\indd {A'} C {B'}$
 \end{enumerate}

The compatibility of the two definitions is shown, e.g., in Section 3  
of \cite{BaldwinShiJapan}.  The following is well known (cf. 3.31 of  
\cite{BaldwinShiJapan}).

\begin{lemm}
\label{cl}
Suppose $A,B$ and $C = A \inter B$ are closed and $\indd A
C B$.  Then $AB$ is closed, i.e. $\bar{AB} = \Abar \cup \Bbar$.
\end{lemm}

The equivalence of $d$-independence and stability theoretic  
independence
was first proved in this generality in \cite{BaldwinShiJapan} but the  
basic setup  comes from \cite{Hrustableplane}.

\jtbnumpar{Fact}
\label{corforkd}
Suppose $T$ satisfies Context~\ref{finalsetup}.
If $C$ is intrinsically closed then for any $A$ and $B$,
$\ind A C B$ if and only if $\indd A C B$.

We give a different proof that is not as involved with the  
intricacies
of amalgamation in the case without finite closures as the one in
\cite{BaldwinShiJapan}.

Suppose for contradiction that $R(A,C,B) \neq \emptyset$.  Then
for $\epsilon$ chosen according to Axiom~\ref{yax},
$\delta(A/B) - \delta(A/BC) > \epsilon$. 
% (by  Fact~\ref{deltadisjoint}).
Now, construct a nonforking sequence $\langle A_i,B_i\rangle$ in
$\tp(AB/C)$.  Since $A$ is not in the algebraic closure of $BC$, no
$A_j$ is in the algebraic closure of the union of $B_i$ for $i<j$.
We will use this fact to show that the types $p_i = \tp(A_i/CB_i)$  
are
$n$-contradictory for some $n$.
If not, for each $n$ there is an $A^*$ which is
common solution for, say $p_1, \dots,   p_n$.
Fix $n$ such that $n\cdot \epsilon >
\delta(A/C)$.
But $\delta(A^*/B_1,
\dots B_n) \leq \delta(A/C) - n \cdot \alpha$ so this implies $A
\subseteq \acl(CB_1 \dots, B_n)$ and this contradiction yields the
result.   The extension property for nonforking types and uniqueness
suffice to deduce the converse from $d$-dependence implies forking
dependence so we finish as in  Lemma 3.35 of \cite{BaldwinShiJapan}.

We extend our notion of dimension to a global real-valued rank on  
types.

\jtbdef Let $p\in S(A)$.
Define $d(p)$ as $d(\abar/A)$ for some (any) $\abar$ realizing $p$.

\jtbdef Let $p_1, p_2 \in S(A)$.
\mbox{}
\begin{enumerate}
\item $p_1$ and $p_2$ are disjoint if for any $\abar_1$, $\abar_2$
realizing
$p_1$, $p_2$,\\ $\icl (A\abar_1) \inter \icl(A\abar_2) \subseteq  
\icl(A)$.
\item $p_1\in S(A)$ and $p_2\in S(B)$
 are disjoint if any pair of nonforking extensions of $p_1$ and $p_2$  
to $AB$ are disjoint.
\end{enumerate}

%In general intrinsic closure does not determine a geometry (exchange
%fails) but we show that on types with dimension $0$, it does.

\begin{lemm} Let $A \subset B$, $p \in S(B)$ and $p|A = q$ and  
suppose
$A$ is intrinsically closed.
 \begin{enumerate}
\item If $d(p) < d(q)$ then $p$ forks over $A$.
\item $q$ is stationary.
%\item If $d(q) = 0$ then $q$ is minimal.
\end{enumerate}
\end{lemm}
\Proof.  1) follows immediately from Fact~\ref{corforkd}; 2) is also
proved in \cite{BaldwinShiJapan} (Lemma 3.38).
%Note that
%a stationary type $q$ is minimal if and only if every forking extension
%of $q$ is algebraic.
%Now,
%for iii),  since $d(q) = 0$, for any
%$\abar$ realizing $q$ in a model $N$, if $\tp(\abar/B)$ forks over $A$,     
%there is
%a finite $C \subseteq \icl_N(A\abar)$ such that $
%C \inter \icl(B) \neq
%\emptyset$.  By Lemma~\ref{iclexc}, $\abar \in \icl(A \cup
%(C \inter \icl(B))$ as required.
\begin{lemm}
\label{orth}
Let $A$ be intrinsically closed, $p_1, p_2\in S(A)$.  If $p_1$ and  
$p_2$ are disjoint and $d(p_1) = 0$ then $p_1$ and $p_2$ are  
orthogonal.
\end{lemm}

\Proof.  If not, there exist sequences
$a_1 \ldots a_k$ and $b_1 \ldots b_m$ of realizations of $p_1$
and $p_2$ respectively,
which are independent over $A$,
 such that
$\dep {\abar} A {\bbar}$.  Since $d(p_1) = 0$,
$d(\abar/A) = 0$ and
$\icl(A\abar) \inter \icl(A\bbar)
\not \subseteq A$.  By Lemma~\ref{cl}, intrinsic closure is a trivial  
dependence relation.
Since the $a_i$ and the $b_j$ are
independent,
this implies that for some $i,j$,
$\icl(A a_i) \inter \icl(A b_j) \not \subseteq A$.  But this  
contradicts
the
disjointness of $p_1$ and $p_2$ and we finish.

The {\em dimensional order property} (DOP) and {\em dimensional
discontinuity property} {DIDIP} are defined in \cite{Shelahbook2nd}.
Either of these conditions implies $T$ has many models in uncountable  
powers.
$T$ has
the eventually non-isolated dimensional order property (eni-dop)
if some type witnessing the dimension order property is not isolated.   
This condition implies that $T$ has the maximal number of countable  
models.
Since $T_{\alpha}$ is not small for irrational $\alpha$, this is not  
new information.  However, the eni-dop seems to be a much more  
intrinsic feature of the construction than the smallness.
(For precise definition see e.g. \cite{Baldwinbook}.)

\begin{thm}
\label{didip}Let $\bK_0$ be a class satisfying  
Context~\ref{finalsetup}.
%and closed under free amalgamation.
Let $T$ be the theory of the generic model for $(\bK_0,\subm)$.  
Suppose further that there is a pair of independent points $B =  
\{x,y\}$ and a nonalgebraic type $p$ with $d(p/B) = 0$ but $d(p/x) >  
0$ and $d(p/y)> 0$.
\begin{enumerate}

\item The theory $T$ has the dimensional order
property.
\item If $p$ is not isolated the theory $T$ has the eni dimensional  
order
property.
\item  The theory $T$ has the  dimensional
discontinuity property.
\end{enumerate}
\end{thm}

\Proof. i)  Let $A = \{a,b\}$ where  $a$ and $b$ are independent  
over
the empty set.
%\abar$ and $\bbar$
%ave the same length and $\abar\bbar$ is an independent sequence.
It
suffices to show that there is a type $p \in S(A)$ with $d(p) = 0$  
and
such that if $\cbar$ realizes $p$, $\dep \cbar {a} b$ and
$\dep \cbar {b} a$.
For then we can construct an
independent sequence of points $a_{i}$ and disjoint copies $p_{i,j}$
over $\{a_i,a_j\}$ which will be pairwise orthogonal by
Lemma~\ref{orth}.
The required type is constructed in Theorem~\ref{mintype}.  ii)  
follows
by the same argument if $p$ is not isolated.

For iii) it suffices to find an independent sequence of sets $B_n$  
for
$n< \omega$ and
$p \in S(B)$ where $B = \union B_n$ such that
$p \orth \union_{n<j} B_n$ for each $j$. Choose $B_n$ and $C_n$ as  
described
at the beginning of the proof of Theorem~\ref{mintype}. Let $B$ be  
the union for $n<\omega$ of $B_n = \{x_n,y_n\}$ with no
relations on $B$.
For each $n$, let $f_n$ map $c_n$ to $c$, $x$ to $x_n$ and $y$ to
$y_n$.  Then $B \union\{c\}$ is as required.
That is,
$d(t(c/B))= 0$ but
$d(t(c/\union_{n<m} B_n)) >0$.

\section{Constructing types of  $d$-rank 0}
\label{findtype}
 We construct  a nonalgebraic type $p$ over a two element set  with  
$d(p)=0$.

\jtbnumpar{Context} We  work with a class $\bK_0$ of finite  
structures as in
Example~\ref{mainexamp}. Thus, $(\bK_0,\subm)$ witnesses
Contex~\ref{finalsetup}. Recall that $\bK$ is the class of models of  
the theory of the generic $M$,
$\Mscr$ is a saturated model of this theory, and $S(\bK)$ is the  
universal class it determines.

Finally, the $\alpha$ parameterizing the dimension function may
 be rational or irrational.  This distinction affects only the  
question
of whether the  type with rank $0$ is isolated  and we discuss that  
when it
arises.

\jtbdef
\label{fullam}
$({\bK_{0}},\subm)$ has the {\em full amalgamation property} if
$B \inter C = A$ and $A \subm B$ imply
$B \bigotimes_{A} C \in {\bK_{0}}$ and $C\subm B \bigotimes_{A} C$.

It is easy to check (Section 4 of \cite{BaldwinShiJapan})that if  
$(\bK_0,\subm)$ is closed
under free amalgamation then it has full amalgamation.

\jtbnumpar{Assumption}
\label{freeamalg}
$({\bK_{0}},\subm)$ has the {\em full amalgamation property}.

\jtbnumpar{Examples} Each of the following classes is closed under  
free amalgamation.
\begin{enumerate}
\item The class $(\bK_{\alpha}, \subm)$ of all finite $L$-structures  
$A$ with $\delta_{1,\alpha}(A)$
hereditarily positive. The resulting theory is $\omega$-stable if  
$\alpha$
is rational and stable if $\alpha$ is irrational.
\item The class yielding the stable $\aleph_0$-categorical  
pseudoplane of
\cite{Hrustableplane}.
\end{enumerate}

The main aim of this section is to establish the following result
which leads easily by Theorem~\ref{didip} to showing the theory of  
the generic model $\Mscr$ has DOP
and DIDIP.

\jtbdef We say $C$ is a {\em primitive} extension of $B$ if $B \subm  
C$ but there is no $B'$ properly between $B$ and $C$ with $B' \subm  
C$.

\begin{thm}
\label{mintype}
%For every $\gamma$ with $0\leq \gamma \leq 1$
There exists a
triple $\{x,y,c\}\in \Mscr$
such that $B = \{x,y\}$ is an independent pair over $\emptyset$ and  
$d(c/xy) = 0$ but
$d(c/x) > 0$, $d(c/y) > 0$ and $c \not \in \acl(x,y)$.
\end{thm}

\Proof.  Fix a discrete structure $B$ with universe $ \{x,y\}$.
We will construct a family $\langle (C_n,x_n,y_n,c_n): n <
\omega\rangle$ of structures in $\bK_0$ which satisfy the following
conditions. Let $B_n =
\{x_n,y_n\}$.  The inequalities in the following discussion
automatically become strict inequalities if $\alpha$ is irrational.

\begin{enumerate}
\item $0 \leq \delta(C_n/B_n) < 1/n$.
\item $(x_n,y_n,c_n)$ is a discrete substructure of $C_n$.
\item $C_n$ is a primitive extension of $B_n$.
%For every $C'$ with $B_n \subseteq C' \subseteq C_n$,
%$\delta(C'/B_n) \geq \delta(C/B_n)$.
\end{enumerate}

Now map each $B_n$ to $B$ and amalgamate the images of the
$C_n$ disjointly over $B$. Then identify
all the $c_n$ as $c$ to form a structure $A$.  Without loss
of generality we can assume $A$ is strongly embedded in $\Mscr$.
Thus, $\icl_{\Mscr}(cB) = A$.
Then $d(c/B) = 0$ but $d(c/x)$ and $d(c/y)$ are both at least one.
Thus $\dep c x {xy}$ and $\dep c y {xy}$.
Since $\delta(C_n/B_n) \geq 0$, for every $n$, $c \not\in \acl(B)$.

\jtbnumpar{Remark} If $\alpha$ is irrational, all the $C_n$ are  
necessary and
 $\tp(c/xy)$ is nonprincipal.
If $\alpha$ is rational, for some $n$, $\delta(C_n/B_n) = 0$. (We  
expand on this
remark after Observation~\ref{obs}.)  The type is principal
but still not algebraic since in this context there are infinitely  
many
copies (in a generic) of a primitive extension with relative  
dimension $0$.
\vskip .2in

The construction of the $C_n$ follows a rather tortured path.  We
first need to consider structures with negative dimension over $B$.

\jtbdef Let $\Ascr =\Ascr_{\alpha}$ be the class of structures of
the form $( A, a, b,e)$
which satisfy the following conditions.
Let $B$ be
the structure with universe $\{a,b\}$ and no relations.
\begin{enumerate}
\item
$A \in \bK_0$.
\item
$\{a,b,e\}$ is the universe of a discrete
substructure of $A$.
\item For each $A'$ with $B \subseteq A'$ and $A'$ properly contained  
in  $A$, $\delta(A')
> \delta(A)$.
\item $-1 < \delta(A/B) \leq 0$.

\end{enumerate}

\jtbnumpar{Observation}
\label{obs}
\begin{enumerate}
\item The
choice of $\delta$ as $\delta_{\alpha}$ makes $\Ascr$ depend on  
$\alpha$.
\item  If the last three conditions are satisfied, the first is as  
well.
%If $\alpha$ is irrational, the inequalities in this definition are
%automatically strict.
\item The last condition implies that $\delta(A/a) > 0$ and  
$\delta(A/b)
> 0$.
\end{enumerate}
\vskip .3in

 We first show
that the set $$X = X_{\alpha}=
\{\beta: \beta = \delta(A/\{a,b\}) \text{ for some } ( A, a,b,e)\in  
\Ascr 
\}$$
is not bounded away from zero.  
If $\alpha$ is irrational,  $0 \not \in X$ so 
$X$ is infinite.  If $\alpha=p/q$ is rational, every element of $X$  
has the form
$(mq-np)/q$ so there cannot be an infinite sequence of members of $X$  
tending
to $0$.  That is, there will be an $A$ with $\delta(A/B) =0$.
As indicated $X$ depends on $\alpha$ (through $\delta =  
\delta_{\alpha}$ and $\Ascr = \Ascr_{\alpha}$.)  But the bulk of the  
proof is uniform in $\alpha$, 
so to enhance readability we keep track of $\alpha$ only for that  
part
 of the proof where the dependence is not uniform.

\jtbnumpar{Construction}
\label{construct}
There are two elementary steps in the construction.  It is easy to  
check
that if the constituent models described here are in $\bK_0$, then so  
is
the result.

\begin{enumerate}
\item If $\delta(A/B) = \beta$
and $\beta \in X$,
and $A^*$ is the free amalgam
over $B$ of $k$ copies of $A$, then $\delta(A^*/B) = k\beta$.

\item Let $(A_1,a_1,b_1,c_1)$ and $(A_2,a_2,b_2,c_2)$ be in $\Ascr$.  

Let $A^*$ be formed
by identifying $b_1$ and $a_2$ and freely amalgamating over that
point. 

\end{enumerate}

\begin{lemm} If $\beta >-1/k$
and $\beta \in X$ then
$k\beta \in X$.
\end{lemm}

\Proof.  Use Construction~\ref{construct} i).

It is straightforward to determine the following properties of the  
second construction.

\begin{lemm}
\label{propconst}
Suppose  $\delta(A_1/\{a_1,b_1\}) = \beta_1$,
$\delta(A_2/\{a_2,b_2\}) = \beta_2$ and $\beta_1, \beta_2 \in X$.   
Let $A^*$ be formed as in Construction~\ref{construct} ii).
\begin{enumerate}
\item $\delta(
A^*/
\{a_1,b_2\})
= \beta_1+\beta_2 +1$.
\item
If
$-2 < \beta_1
+\beta_2 \leq -1$
then $\beta_1 +\beta_2 +1 \in X$ and
$\langle A^*,a_1,b_2,c_1\rangle \in \Ascr$.
\item
If
$-1 \leq \beta_1
+\beta_2< -1+1/n$ then
\begin{enumerate}
\item
$0 \leq
\delta(A^*/\{a_1,b_2\})< 1/n$.
\item
$\delta(A^*/a_1) \geq 1$
and
$\delta(A^*/b_2) \geq 1$.
\end{enumerate}
\end{enumerate}
\end{lemm}
\proof.
The key observations for 1)and thus 2) and 3a) is that for any $B  
\subseteq A_1 \subseteq A^*$,
$$\delta(A'/\{a_1,b_2\})=\delta(A'\inter A_1/\{a_1,b_1\})  
+\delta(A'\inter A_2/\{a_2,b_2\}) + 1.$$

For 3b) we need the further remark:
$$\delta(A'/a_1)=\delta(A'/b_2)=\delta(A'/\{a_1,b_2\}) +1.$$

\begin{lemm}
\label{stageone}
If $L$ contains a single binary relation and $\bK_0 = \bK_{\alpha}$,
then $X$ is not empty.
\end{lemm}

%The construction depends on $\alpha$ so for the duration of the  proof
%we refine our notation to write $$Y_{\alpha}=\{A:
%\Ascr_{\alpha}A \in \Ascr \ \&\  -1< \delta_{\alpha}(A/B)\leq  0\}.$$

\proof.
It suffices to show that each $\Ascr_{\alpha}$ is nonempty for $0<  
\alpha\leq
1$.  The construction is somewhat ad hoc and proceeds by a number of
cases depending on $\alpha$.  Thus to establish Lemma~\ref{stageone}  
we will use
the notations $\Ascr_{\alpha}, \delta_{\alpha}$.
These constructions are very specific to graphs.  The second author  
has
an alternative argument which avoids the dependence on $\alpha$.   
However, it passes through hypergraphs and has it own computational  
complexities.

\jtbnumpar{Case 1} $3/4 < \alpha < 1$:  Let $A_1$ be the structure
obtained by adding to $\{a,b,e\}$ two points $b_1,b_2$ such that
$b_1$ is connected to $a$ and $e$ while $b_2$ is connected to
$b$ and $e$.  Then

$$-1 < \delta_{\alpha}(A_1/B) = 3 - 4\alpha < 0$$
for the indicated $\alpha$ and $(A_1,a,b,e) \in \Ascr_{\alpha}$.

\jtbnumpar{Case 2} $2/3\leq \alpha < 4/5$:  Let $A_2$ be the  
structure
obtained by adding to $\{a,b,e\}$ two points $b_1,b_2$ such that
$b_1$ is connected to $a$, $b$, and $e$ while $b_2$ is connected to
$b$ and $e$.  Then

$$-1 < \delta_{\alpha}(A_2/B) = 3 - 5\alpha < 0$$
for the indicated $\alpha$ and $(A_2,a,b,e) \in \Ascr_{\alpha}$.

\jtbnumpar{Case 3} $0< \alpha < 2/3$:  Let $A_{n,k}$ be the  
structure
obtained by adding to $\{a,b,e\}$ both
$n$ points $a_1,\ldots,a_n$ such
that each $a_i$ is connected to $a$, $b$, and $e$ and
$k$ points $b_1,\ldots,b_k$ such
that each $b_i$ is connected to all the $a_i$.

Then $\delta_{\alpha}(A_{n,k}/B) =
n +k +1 - (nk + 3n)\alpha$.  We say $\alpha$ is {\em acceptable} for  
$n$ and $k$ if the following inequality is satisfied.
$$\ell_{n,k} =
{{n +k +1}\over {nk + 3n}}
< \alpha <
{{n +k +2}\over {nk + 3n}}
= u_{n,k}.$$

To show that if $\alpha$ is acceptable for $n$ and $k$, then  
$(A_{n,k},a,b,e) \in \Ascr_{\alpha}$ we need several claims.

\jtbnumpar{Claim 1}  For each $k$,
\begin{enumerate}
\item
$u_{n+1,k} > \ell_{n,k}$,
\item $\ell_{n+1,k} < \ell_{n,k}$,
\item
$\lim_{n\rightarrow \infty} \ell_{n,k} = 1/(k+3)$.
\end{enumerate}

Claim 1 is established by routine computations.

\jtbnumpar{Claim 2}  For every $\alpha$ that is acceptable for $n$  
and $k$,
if $B \subseteq A'\subseteq A_{n,k}$,
$\delta_{\alpha}(A'/B) \geq \delta_{\alpha}(A_{n,k}/B)$.

To see this, note that any such $A'$, for  some $m \leq n$ and
$\ell \leq k$, either $A'$ has the form $A_{m,\ell}$ or
the form $B_{m,\ell}$, where $B_{m,\ell}$ is the structure obtained  
by omitting the
element $e$ from $A_{m,\ell}$.  Now note that if  
$\delta_{\alpha}(B_{m,\ell}/B) < 0$ then
$\delta_{\alpha}(B_{m,\ell}/B) \geq \delta_{\alpha}(B_{m+1,\ell}/B) $  
and
$\delta_{\alpha}(B_{m,\ell}/B) \geq\delta_{\alpha}(B_{m,\ell+1}/B) $.   
The same assertion
holds when $A_{m,\ell}$ is substituted for $B_{m,\ell}$.  Finally,
$\delta_{\alpha}(B_{n,k}/B)\geq \delta_{\alpha}(A_{n,k}/B) $.  These  
three observations
yield the second claim.

\relax From these two claims we see that for each $\alpha$, there is
a pair $n,k$ with $A_{n,k}\in \Ascr_{\alpha}$.
The remainder
of the argument does not depend on $\alpha$ so we return to the use
of the notation $X$ and $\Ascr$.

\begin{lemm}
\label{smalld}
For every $n$ there is an element $\beta$ of $X$ with $\beta>-1/n$.
\end{lemm}

\Proof.  If not, fix the least $n$ such that all elements of $X$ are
at most\\
$-1/(n+1)$ and fix $\beta_0\in X$ with $-1/n < \beta_0 \leq  
-1/(n+1)$.  
(If $\beta_0 = -1/(n+1)$, $\beta_1 = 0$ and we finish.)  Define by
induction $\beta_{\ell +1} = (n+1) \beta_{\ell} + 1$.  Combining the  
two
elementary steps we see that each $\beta_{\ell} \in X$.  Let $
\beta'_{\ell}$ be the distance between $-1/n$ and $\beta_{\ell}$.   
That is,
$\beta'_{\ell}=|-1/n -\beta_{\ell}|=
1/n + \beta_{\ell}$.
Now $\beta_{\ell} \leq -1/(n+1)$ if and only if
$\beta'_{\ell} \leq 1/(n)(n+1)$.

But $$
\beta'_{\ell+1} =
1/n +(n+1) \beta_{\ell}+ 1 = (n+1) \beta'_{\ell}.$$
So $$
\beta'_{\ell} =
(n+1)^{\ell} \beta'_{0}.$$
As $\beta'_0 > 0$, for sufficiently large $\ell$,
$\beta'_{\ell} > 1/(n)(n+1)$ so $\beta_{\ell} > -1/(n+1)$ as  
required.

With a few more applications of our fundamental constructions, we  
can
find the $C_n$ needed for Theorem~\ref{mintype}.

By applying Construction~\ref{construct} i) and Lemma~\ref{smalld}
for any $n$, and $i = 1,2$ we
can find
$(A^n_1,x^n_1,y^n_1,c^n_1)$ and $(A^n_2,x^n_2,y^n_n,c^n_2)$
containing $B^n_i = \{x^n_i,y^n_i\}$
such that $\{x^n_i,y^n_i,c^n_i\}$ is discrete and  
$\delta(A^n_i/B^n_i) =
\beta^n_i$ with
$ - 1 < \beta^n_1 + \beta^n_2 < -1+ 1/n$.

To construct $A^n_1$,
choose using
Lemma~\ref{smalld} a $(D^n,x^n_1,y^n_1, c^n_1) \in \Ascr$ with $-1/n  
<
\delta(D^n/B^n_1)
\leq 0$.  Take an appropriate number, $k$, of
copies of $D^n$ over $B^n_1$ and apply
Construction~\ref{construct} i) to form $A^n_1$
 with
$$-1 < k\delta(D^n/B^n_1) = \delta(A^n_1/B^n_1)=\beta^n_1 < -1 +  
1/n$$
 and choose $c^n_1 \in A^n_1$ so that
$(x^n_1,y^n_1,c^n_1)$ is discrete.
By Lemma~\ref{smalld} again choose $(A^n_2,x^n_2,y^n_2, c^n_2)\in  
\Ascr$
with $$-(\beta^n_1 +1)/2 < \delta(A^n_2/B^n_2) = \beta^n_2 < 0.$$
Now apply
Construction~\ref{construct} ii) to
$(A^n_1,x^n_1,y^n_1,c^n_1)$ and $(A^n_2,x^n_2,y^n_n,c^n_2)$
to form $(C_n,x_n,y_n,c_n)$
where $x_n = x^n_1$, $y_n = y^n_2$, and $c_n = c^n_1$.  Denote
$\{x_n,y_n\}$ by $B_n$.  Then
$0<\delta(C_n/B_n) =
1 +\beta^n_1 +\beta^n_2 < 1/n$.  Each $C_n$ contains a discrete set
$\{x_n,y_n,c_n\}$ and the third property of the $C_n$ follows using
the second part of Lemma~\ref{propconst}.
This completes the construction of the type of $d$-rank
0.

Using the argument for constructing $A^n_1$,
 we easily
show the following density result.

\begin{cor}
\label{dense}
For any  $\gamma, \delta$ with $-1 \leq \gamma < \delta < 0$ there is  
a $(D,a,b,e) \in
\Ascr$ with $\gamma < \delta(D/\{a,b\}) < \delta$.
\end{cor}

 The restriction to one-types in the following lemma is solely for  
ease of
presentation.

\begin{lemm}  Suppose $A\subseteq M \sat T_{\alpha}$ is  
intrinsically
closed and $p_1,p_2
\in S_1(A)$ are disjoint. If $0<d(p_i)$ for $i = 1,2$ then $p_1  
\nperp p_2$.
\end{lemm}

\Proof.  Clearly if $p_1$ and $p_2$ are not disjoint or if
there is an edge between realizations of the two types, they
are not orthogonal. Let $a_1,a_2$ realize $p_1,p_2$
and suppose for contradiction that $p_1$ and $p_2$ are orthogonal
and $d(a_1a_2/A) = d(a_1/A) + d(a_2/A)
= \beta> 0$. In particular, there is
no edge linking
$a_1$ and $a_2$.   By Lemma 3.25 of \cite{BaldwinShiJapan} there are  
finite $A_1
\supseteq a_1a_2$ and $A_0 \subset A$ with $\beta \leq \gamma =  
\delta(A_1/A_0)<
\beta +1$.  Lemma~\ref{dense} allows us to choose a finite $B  
\supseteq
\{a_1,a_2\}$ with $$-1 < \delta(B/\{a_1,a_2\}) <\beta - \gamma <0.$$

Then $Ba_1a_2$ is in $\bK_0$.  By full amalgamation we can freely  
amalgamate $B$ with $AA_1$ over
$\{a_1,a_2\}$ inside $\Mscr$.
Then $d(a_1a_2/A) \leq \delta(A_1B/A_0)$.   Note
$\delta(B/A_1A_0) = \delta(B/\{a_1,a_2\}) < \beta-\gamma$.
So $$\delta(A_1B/A_0) = \delta(B/A_1A_0) + \delta(A_1/A_0) <  
\beta.$$
This contradicts $d(a_1a_2/A)= \beta$ so we conclude $p_1 \nperp  
p_2$.

Using the Lemmas~\ref{orth} and \ref{dense} it is easy to see

\begin{cor}  In $T_{\alpha}$, 
\begin{enumerate}
\item For disjoint $p_1, p_2$, $p_1 \perp p_2$ if and only if $d(p_1)  
= 0$ or $d(p_2) = 0$.
\item Every regular type satisfies $d(p) = 0$.
\end{enumerate}
\end{cor}

Our construction yields some further information.

\jtbdef   The type $p \in S(A)$ is {\em minimal} if $p$ is not  
algebraic but
for any formula $\phi(x,\bbar)$ either
$p \cup \{\phi(x,\bbar)\}$
or
$p \cup \{\neg\phi(x,\bbar)\}$
is algebraic.

\jtbdef The type $p \in S(A)$ is $i$-minimal if for every $\abar$
realizing $p$, if $c \in \icl(A\abar)$, $\icl(Ac) = \icl(A\abar)$.

\begin{thm}  If $p$ is constructed as in Lemma~\ref{mintype} then $p$  
is minimal and trivial.
\end{thm}

\proof.  If $d(p) = 0$ and $p$ is i-minimal then $p$ is
minimal.  We constructed $p$ so that $d(p) = 0$ but the fact that  
each $C_n$ is primitive over $B$ and $A$ is intrinsically closed  
guarantees that $p$ is $i$-minimal and we finish. 

Clearly, $d(p) = 0$ does not imply $p$ is minimal.  For,
if $d(a/A) = d(b/A) = 0$ then $d(ab/A) = 0$ but if, for example, $a$  
and
$b$ are independent $\tp(ab/A)$ is not minimal.

\section{The Finite Cover Property}
\label{fcp}
In this section we show that for classes as described in  
Example~\ref{mainexamp} with the full amalgamation
property, and in particular for $(\bK_{\alpha},\subm)$, the theory  
of
the generic does not have the finite cover property.  We rely on the
following characterization due to Shelah  
\cite[II.2.4]{Shelahbook2nd}.

\jtbnumpar{Fact}
\label{fcpfact}
If $T$ is a stable theory with the finite cover
property then there is a formula $\phi(\xbar,\ybar,\zbar)$ such that
\begin{enumerate}
\item For every $\cbar$, $\phi(\cbar,\ybar,\zbar)$ defines an
equivalence relation.
We call this relation $\cbar$-equivalence.
\item For arbitrarily large $n$, there exists $\cbar_n$ such that  
the
equivalence relation defined by $\phi(\cbar_n,\ybar,\zbar)$ has  
exactly
$n$ equivalence classes.
\end{enumerate}

Here is some necessary notation.

\jtbdef Let $A,B$ be finite substructures of $M$ with $A \subseteq  
B$
then
\begin{enumerate}
\item
$\chi_M(B/A)$ is the number of distinct copies of $B$ over $A$
  in $M$.
\item $\chi^*_M(B/A)$ is the supremum of the cardinalities of  
maximal
  families of disjoint (over $A$) copies of $B$ over $A$ in $M$.
\end{enumerate}

\jtbdef $(A,B)$ is a {\em minimal pair} if $\delta(B/A) < 0$ and for  
every 
$B'$, with $A \subseteq B' \subseteq B$, $\delta(B/A) <  
\delta(B'/A)$.

The next result is proved in \cite{BaldwinShiJapan}.

\jtbnumpar{Fact}
\label{intrinsfew} There is a function $t$ taking pairs of integers  
to
integers such that
if $A\subi B$ then 
 for any $N\in \bK$ and any embedding
$f$ of $A$ into $N$,
$\chi_N(fB/fA)\leq t(|A|,|B|)$.

There is an easy partial converse to this result.

\begin{lemm}
\label{manystrong}

For any $M \in \bK_0$, if
$\chi^*_M(B/A) > t(|A|,|B|)$ then
$A \subm B$.

\end{lemm}

\proof.  Suppose some $B'$ with $A \subseteq B$ satisfies $A \subi  
B'$.
Then there are more than $t(|A|,|B|)$ disjoint copies of $B'$ over  
$A$ in $M$ contradicting Fact~\ref{intrinsfew}.

We also need the finer analysis of the intrinsic closure carried out
in \cite{BaldwinShelahran}.  In fact, this argument depends on the  
slightly finer notion of a {\em semigeneric} which is defined in  
\cite{BaldwinShelahran}.  The crucial facts from  
\cite{BaldwinShiJapan} and \cite{BaldwinShelahran} are the following.

\jtbnumpar{Fact}
\label{locclosure} If $(\bK_0,\subm)$ satisfies  
Context~\ref{finalsetup} and has the
full amalgamation property then the theory of the generic $T$  
satisfies
\begin{enumerate}
\item All models of $T$ are semigeneric.
\item $T$ is stable.
For any formula $\phi(x_1 \ldots x_r)$ there is an
integer $\ell = \ell_{\phi}$, such that for any semigeneric $M\in  
\bK$
and any $r$-tuples $\abar$
and $\abar'$ from $M$ if
$\icl^{\ell_{\phi}}_M(\abar)
\iso
\icl^{\ell_{\phi}}_M(\abar')$ then
$M\sat \phi(\abar)$
if and only if
$M\sat \phi(\abar')$.
\end{enumerate}

\begin{thm}  Let the language $L$ contain only binary relation  
symbols.
If $(\bK_0,\subm)$ satisfies Context~\ref{finalsetup} and has the
full amalgamation property then the theory of the generic $T$ does  
not
have the finite cover property.
\end{thm}

\Proof.  Suppose not.  We know $T$ is stable so there is a formula  
$\phi$ satisfying
the conditions of Fact~\ref{fcpfact}.  Each model of $T$ is  
semigeneric.
Choose $\ell = \ell_{\phi}$ as in Fact~\ref{locclosure} so that the
isomorphism type of
$\icl^{\ell}_M(\cbar,\abar,\bbar)$
determines the
truth of $\phi(\cbar,\abar,\bbar)$ for any triple of  
$\cbar,\abar,\bbar$
of appropriate length.  For any $n$ choose $m$ sufficiently large  
with
respect to the maximal cardinality of
$\icl^{\ell}_M(\cbar,\abar,\bbar)$ and $n$ so
that applying the pigeonhole principle  and Ramsey's theorem
we can choose $\cbar_m$ so that the
$\langle \abar_i:i< n\rangle$ are pairwise $\cbar_m$-inequivalent and  
for $i < n$
letting $A_i = \icl_M^\ell(\cbar,\abar_i)$ and $C =  
\icl^{\ell}_M(\cbar)$
the following property $P(C)$ holds.
\begin{enumerate}
\item for all $i,j$, $A_i \iso_{C} A_j$
\item for $i<j$, $A_0A_1 \iso_{C} A_iA_j$.
 \end{enumerate}

If $n>t(k,|A_0|)!$ for $k< |A_0|$, applying  the $\Delta$-system  
Lemma we can find ${\hat C}$ with $C \subseteq {\hat C} \subseteq  
A_0$ such that (without loss of generality) the $A_i$ are disjoint  
over ${\hat C}$.  By appropriate choice of $n$, depending only on  
$|A_0|$, $|C|$, we may assume
that $p({\hat C})$ holds. 
By Fact~\ref{manystrong}, ${\hat C} \subm A_0$.  We claim in fact  
that the
structure imposed on $A_0A_1$ is $A_0 \otimes_{\hat C} A_1$.  If  
not,
$R(A_0,{\hat C},A_1)$ is nonempty.  Let $E_i$ denote the substructure   
of M with
universe $\bigcup _{j<i} A_j$.  By Axiom~\ref{yax}~iii) for
sufficiently large $k$, $\delta(A_k/E_k) < 0$.  There is a minimal  
pair
$(E'_k,A'_k)$ with $E'_k \subseteq E_k$ and $A'_k \subseteq A_k$.
But then for each $j > k$ there is a copy $A'_j$ of $A'_k$, contained  
in
$A_j$ and isomorphic to $A'_k$ over $E_k$ (since the language is  
binary).  This contradicts the bound
on the number of copies of a minimal pair, Fact~\ref{intrinsfew}.
Thus we establish the claim.  But now we have $E_{i+1} \iso E_i
\otimes_C A_0$.  Since $(\bK_0,\subm)$ has full amalgamation, this
construction can be carried on indefinitely.  But the definition of
$\ell_{\phi}$ guarantees that the $\abar_i$ represent distinct
$\cbar$-equivalence classes and this contradicts the hypothesis that
there are only finitely many $\cbar$-equivalence classes.

\jtbnumpar{Conclusion}  The arguments in the paper are fully worked  
out only for languages with binary relation symbols.  For Section 4,  
this is just a matter of easing notation; slight modifications of the  
argument work for any finite relational language.  The combinatorial  
arguments in Section 3 are sufficiently complicated that the proof is  
the general case is less clear.  But it would be quite surprising if  
the restriction to a binary language is actually necessary.

\bibliography{ssgroups}
\bibliographystyle{plain}
\end{document}

%% file: jbmac.tex
\makeatletter
\def\@begintheorem#1#2{\it \trivlist \item[\hskip
\labelsep{\bf #2\ #1.}]}
\def\@opargbegintheorem#1#2#3{\it \trivlist
      \item[\hskip \labelsep{\bf #1\ #2\ (#3).}]}
\makeatother

\newtheorem{them}{Theorem}[section]

\newtheorem{lemm}[them]{Lemma}
\newtheorem{cor}[them]{Corollary}

\newtheorem{thm}[them]{Theorem}

\newcommand{\jtbnumpar}[1]{\refstepcounter{them}
\trivlist
\item[\hskip \labelsep{\bf \thethm \ #1.}]}

\newcommand{\jtbdef}{\jtbnumpar{Definition}}
\def\jtbnot{\jtbnumpar{Notation}}
\def\PROOF #1.{\par\noindent{\it Proof#1}.\ \ignorespaces}
\def\proof #1.{\par\noindent{\it Proof#1}.\ \ignorespaces}

\def\subm{\leq}

\def\Mscr{{\cal M}}

\def\Ascr{{\cal A}}

\let\?=\joinrel

\let\leftv=^
\let\rightv=^

\def\k/{\kern.2em}    % for eg Theorem 1.2\k/ii)
\def\acl{\mathop{\rm acl}}

\def\acl{{\rm acl}}     % Ch 17
\def\down{\smash{\mathchar"0223}}

\let\union=\cup             %

\edef\bigcup{\mathop{\textstyle\mathchar\the\bigcup}}
\let\bigunion =\bigcup

\let\inter=\cap             %

\edef\bigcap{\mathop{\textstyle\mathchar\the\bigcap}}

\edef\bigwedge{\mathop{\textstyle\mathchar\the\bigwedge}}

\edef\bigvee{\mathop{\textstyle\mathchar\the\bigvee}}

\edef\sum{\mathop{\textstyle\mathchar\the\sum}}
\def\ind #1#2#3{#1 \mathbin{\mathop{\down}_{#2}} #3}
\def\indd #1#2#3{#1 \mathbin{\mathop{\down^d}_{#2}} #3}
\def\dep #1#2#3{#1 \mathbin{\mathop{\not\mathrel{\mkern-2.5mu}\down}_{#2}} #3}
\let\orth=\dashv

\def\nperp {\not\perp}        %
\def\math&{\ \& \ }

\def\force {\mathrel^\joinrel\rightarrow}
\def\force {\mathrel{\scriptstyle\mathrel^\joinrel\rightarrow}}
\def\forceq {\mathrel{\mathop{\force}\limits_{\textstyle\texsim}}}
\def\forceq{\mathrel^\joinrel
 \mathrel{\mathop{\rightarrow}\limits_{\smash{\textstyle\texsim}}}}
\def\forceq{\mathrel{\scriptstyle\mathrel^\joinrel
 \mathrel{\smash{\mathop{\rightarrow}\limits_{\smash{\raise
 2pt\hbox{$\scriptstyle\texsim$}}}}}}}

\let\exclaim=!                 %
\let\iso\approx             %
\let\texsim=\sim         %
\let\conj\sim             %
\def\conjp #1 {\conj_{#1}}     %
\let\sim\simeq             %
\let\neg=\lnot             %

\def\0bar{\bar 0}         %
\def\1bar{\bar 1}         %
\def\abar{\overline a}

\let\sat=\models                %

\def\bbar{\overline b}
\def\cbar{\overline c}

\def\hbar{\overline h}

\def\xbar{\overline x}
\def\ybar{\overline y}
\def\zbar{\overline z}

\def\Abar{\bar A}
\def\Bbar{\bar B}
\def\Cbar{\bar C}

\def\tp{\rm tp}

\def\Proof{Proof}

\def\th{\rm Th}

\def\bigunion{\union}

\def\text#1{\ifmmode\leavevmode\hbox{#1}\else
   \typeout{Warning: \string\text \space used outside math mode!}
   \begingroup\hbox{#1}\endgroup\fi}

%% file: 567.bbl
\begin{thebibliography}{1}

\bibitem{Baldwinbook}
J.T. Baldwin.
\newblock {\em Fundamentals of Stability Theory}.
\newblock Springer-Verlag, 1988.

\bibitem{BaldwinShelahran}
J.T. Baldwin and S.~Shelah.
\newblock Randomness and genericity.
\newblock in preparation.

\bibitem{BaldwinShiJapan}
J.T. Baldwin and Niandong Shi.
\newblock Stable generic structures.
\newblock 1993.
\newblock submitted.

\bibitem{Hrustableplane}
E.~Hrushovski.
\newblock A stable $\aleph_0$-categorical pseudoplane.
\newblock preprint, 1988.

\bibitem{Shelahbook2nd}
S.~Shelah.
\newblock {\em Classification {T}heory and the {N}umber of {N}onisomorphic
  {M}odels}.
\newblock North-Holland, 1991.
\newblock second edition.

\bibitem{ShelahSpencer}
S.~Shelah and J.~Spencer.
\newblock Zero-one laws for sparse random graphs.
\newblock {\em Journal of A.M.S.}, 1:97--115, 1988.

\end{thebibliography}
